\documentclass{amsproc}

\usepackage{amssymb}

\usepackage{graphicx}

\usepackage[cmtip,all]{xy}

\usepackage{}
\usepackage{amsmath}
\usepackage{amsfonts}
\usepackage{amssymb}

\usepackage{amssymb,amsmath,amsxtra,graphicx}
\usepackage{times}
\usepackage{bbold}
\usepackage{stmaryrd}
\usepackage{array}
\usepackage{dsfont}
\usepackage{url}
\usepackage{tikz}
\usepackage[T1]{fontenc}

\newtheorem{theorem}{Theorem}

\theoremstyle{plain}

\newtheorem*{thmA}{Theorem A}

\newtheorem*{conjA}{Conjecture A}

\newtheorem*{conjB}{Conjecture B}

\newtheorem{lemma}[theorem]{Lemma}

\newtheorem{proposition}[theorem]{Proposition}

\newtheorem*{thm}{Theorem MZ}

\numberwithin{theorem}{section}

\begin{document}
\title[Gaussian Riemann Derivatives]{Gaussian Riemann Derivatives}
\author{J. Marshall Ash}
\address{Department of Mathematics, DePaul University\\Chicago, IL 60614}
\email{mash@depaul.edu}
\author{Stefan Catoiu}
\address{Department of Mathematics, DePaul University\\Chicago, IL 60614}
\email{scatoiu@depaul.edu}
\author{Hajrudin Fejzi\'{c}}
\address{Department of Mathematics, California State University, San Bernardino, CA 92407}
\email{hfejzic@csusb.edu}
\thanks{September 26, 2022. This paper is in final form and no
version of it will be submitted for publication elsewhere.}
\date{September 26, 2022}

\begin{abstract}
J. Marcinkiewicz and A. Zygmund proved in 1936 that, for all functions $f$ and points $x$, the existence of the $n$th Peano derivative $f_{(n)}(x)$ is equivalent to the existence of both $f_{(n-1)}(x)$ and the $n$th generalized Riemann derivative $\widetilde{D}_nf(x)$, based at $x,x+h,x+2h,x+2^2h,\ldots ,x+2^{n-1}h$.
For $q\neq 0,\pm 1$, we introduce: two $q$-analogues of the $n$-th Riemann derivative ${D}_nf(x)$ of~$f$ at~$x$, the $n$-th Gaussian Riemann derivatives ${_q}{D}_nf(x)$ and ${_q}{\bar D}_nf(x)$ are the $n$-th generalized Riemann derivatives based at $x,x+h,x+qh,x+q^2h,\ldots ,x+q^{n-1}h$ and $x+h,x+qh,x+q^2h,\ldots ,x+q^{n}h$; and one analog of the $n$-th symmetric Riemann derivative ${D}_n^sf(x)$, the $n$-th symmetric Gaussian Riemann derivative ${_q}{D}_n^sf(x)$ is the $n$-th generalized Riemann derivative based at $(x),x\pm h,x\pm qh,x\pm q^2h,\ldots ,x\pm q^{m-1}h$, where $m=\lfloor (n+1)/2\rfloor $ and~$(x)$ means that $x$ is taken only for $n$ even.
We provide the exact expressions for their associated differences in terms of Gaussian binomial coefficients; we show that the two $n$th Gaussian derivatives satisfy the above classical theorem, and that the $n$th symmetric Gaussian derivative satisfies a symmetric version of the theorem; and we conjecture that these two results are false for every larger classes of generalized Riemann derivatives, thereby extending two recent conjectures by Ash and Catoiu, both of which we update by answering them in a few cases.
\end{abstract}
\maketitle

\noindent
The first family of generalized derivatives was invented by Riemann in the mid 1800s;~see~\cite{R}. For a~positive integer $n$, the $n$-th \emph{Riemann derivative} of a function $f$ at $x$~is~defined~by~the~limit
{\small\[
D_nf(x)=\lim_{h\rightarrow 0}\frac 1{h^n}\sum_{k=0}^n(-1)^k{n\choose k}f(x+(n-k)h).
\]}The above sum is denoted by $\Delta_n(x,h;f)$ and called the $n$-th \emph{Riemann difference} of $f$ at $x$ and $h$. The sequence of Riemann differences satisfies the recursive relations:
\[\small
\begin{aligned}
\Delta_1(x,h;f)&=f(x+h)-f(x),\\
\Delta_n(x,h;f)&=\Delta_{n-1}(x+h,h;f)-\Delta_{n-1}(x,h;f), \qquad (n\geq 2).
\end{aligned}
\]

The Riemann derivatives were generalized by Denjoy in \cite{D} (1935). An $n$-th \emph{generalized Riemann derivative} of a function $f$ at $x$ is defined by the following~limit:
{\small\[
D_{\mathcal{A}}f(x)=\lim_{h\rightarrow 0}\frac 1{h^n}\sum_{k=0}^{\ell }A_kf(x+a_kh).
\]}The above sum is denoted by $\Delta_{\mathcal{A}}(x,h;f)$ and called an $n$-th \emph{generalized Riemann difference}. Its data vector $\mathcal{A}=\{A_0,\ldots ,A_{\ell };a_0,\ldots ,a_{\ell }\}$ for which the $A_k$ are non-zero and the $a_k$ are distinct is required to satisfy the $n$-th \emph{Vandermonde relations} {\small $\sum_{k=0}^{\ell}A_ka_k^j=\delta_{j,n}\cdot n!$}, for $j=0,1,\ldots ,n$. The points $x+a_0h,\ldots ,x+a_{\ell }h$ are called \emph{base points} of either the derivative or the difference.
The Vandermonde linear system is consistent precisely when $\ell \geq n$, and has a unique solution $(A_0,\ldots ,A_{\ell })$ precisely when $\ell =n$. In this case, the generalized Riemann derivative is considered to be \emph{without~excess}. Unless otherwise stated, all generalized Riemann derivatives in this paper will be without excess.

More examples of generalized Riemann derivatives include the $n$-th \emph{symmetric Riemann derivative} $D_n^sf(x)$ whose associated $n$-th \emph{symmetric Riemann difference},
{\small\[
\Delta_n^s(x,h;f)=\sum_{k=0}^n(-1)^k{\binom nk}f\left(x+\left(\frac n2-k\right)h\right),
\]}
satisfies the recursive relations:
\[\small
\begin{aligned}
\Delta_1^s(x,h;f)&=f(x+h/2)-f(x-h/2),\\
\Delta_n^s(x,h;f)&=\Delta_{n-1}^s(x+h/2,h;f)-\Delta_{n-1}^s(x-h/2,h;f),\qquad (n\geq 2).
\end{aligned}
\]
In general, an $n$-th generalized Riemann difference {\small $\Delta_{\mathcal{A}}(x,h;f)$} is \emph{symmetric} if it satisfies {\small $\Delta_{\mathcal{A}}(x,-h;f)=(-1)^n\Delta_{\mathcal{A}}(x,h;f)$}, and is \emph{even} or \emph{odd} if it is symmetric of even or odd order~$n$. The symmetry of a $n$-th generalized Riemann derivative without excess is equivalent to the symmetry of the set $\{a_0,\ldots ,a_n\}$ relative to the origin.
More generalized Riemann differences are obtained by scaling. A \emph{scale} by $r$ of an $n$-th generalized Riemann difference $\Delta_{\mathcal{A}}(x,h;f)$, where $\mathcal{A}=\{A_k;a_k\}$, is the $n$-th generalized Riemann difference $\Delta_{\mathcal{B}}(x,h;f)$, where $\mathcal{B}=\{B_k=r^{-n}A_k;b_k=ra_k\}$. Moreover, $f$ is $\mathcal{A}$-differentiable at $x$ if and only if $f$ is $\mathcal{B}$-differentiable at $x$ and $D_{\mathcal{A}}f(x)=D_{\mathcal{B}}f(x)$.

The second family of generalized derivatives was introduced by Peano in~\cite{P}~(1892) and further developed by de la Vall\'ee Poussin in \cite{dlVP} (1908). A function $f$ has $n$ \emph{Peano derivatives} at $x$ if there exist numbers {\small $f_{(0)}(x):=f(x),f_{(1)}(x),\ldots f_{(n)}(x)$} such that
{\small\[
f(x+h)=f_{(0)}(x)+f_{(1)}(x)h+f_{(2)}(x)\frac {h^2}{2!}+\cdots +f_{(n)}(x)\frac {h^n}{n!}+o(h^n),
\]}as $h$ approaches zero. The number $f_{(n)}(x)$ is the $n$-th \emph{Peano derivative} of $f$ at $x$. Its existence assumes the existence of every lower order Peano derivative of $f$ at $x$.

A function~$f$ is said to have $n$ \emph{symmetric Peano derivatives} at $x$ if there exist real numbers {\small $f_{(0)}^s(x),f_{(1)}^s(x),\ldots ,f_{(n)}^s(x)$} such that $f_{(0)}^s(x)=f(x)$, for $n$ even, and
{\small\[
\frac 12\{f(x+h)+(-1)^nf(x-h)\}=f_{(0)}^s(x)+f_{(1)}^s(x)h+\cdots +f_{(n)}^s(x)\frac {h^n}{n!}+o(h^n),
\]}as $h$ approaches zero. The number {\small $f_{(n)}^s(x)$} is the $n$-th \emph{symmetric Peano derivative} of $f$ at~$x$. Replacing $h$ with $-h$ in the above displayed equation yields {\small $f_{(n-1)}^s(x)=f_{(n-3)}^s(x)=f_{(n-5)}^s(x)=\cdots =0$}. In this way, if $f$ has $n$ symmetric Peano derivatives at $x$, then $f$ has symmetric Peano derivatives at $x$ of orders $n-2, n-4$, and so on.

The following three simple facts about Peano and generalized Riemann derivatives are useful to remember:
\begin{itemize}
\item The existence of the $n$-th Peano derivative $f_{(n)}(x)$ implies the existence of every $n$-th generalized Riemann derivative $D_{\mathcal{A}}f(x)$ and $f_{(n)}(x)=D_{\mathcal{A}}f(x)$.
\item The existence of an $n$-th generalized Riemann derivative does not enjoy the nice property of forcing the existence of lower order derivatives.
\item The existence of an $n$-th generalized Riemann derivative at $x$ is not enough to guarantee that every $n-1$ times Peano differentiable function at $x$ is also $n$ times Peano differentiable at~$x$.
\end{itemize}
We illustrate the last bullet property by looking at the function $g(x)=x^n\,\mbox{sgn}\,x$ at $x=0$. Since $g(h)=o(h^{n-1})$, $g_{(i)}(0)=0$, for $0\leq i\leq n-1$. Since the parity of~$g$ is opposite to the parity of $n$, $D_n^sg(0)=0$. Finally, $g_{(n)}(0)=\lim_{h\rightarrow 0}R(h)$, where $R(h)=n!{g(h)}/{h^n}$, does not exist since $\lim_{h\rightarrow 0^+}R(h)=n!$ but $\lim_{h\rightarrow 0^-}R(h)=-n!$.

\subsection*{Motivation} Our main motivation comes from the following theorem of classical real analysis from 1936, due to Marcinkiewicz and Zygmund in \cite{MZ}, which in \cite{ACF} was called Theorem~MZ, and which is itself motivated here by the above third bullet property and its associated example:

\begin{thm}[\cite{MZ}, $\S 10$, Lemma~1] Let $\widetilde{D}_n$ be the $n$-th generalized Riemann derivative based at $x,x+h,x+2h,x+4h,\ldots ,x+2^{n-1}h$. Then, for each function $f$ and point $x$,
\[
\text{both $f_{(n-1)}(x)$ and $\widetilde{D}_nf(x)$ exist $\Longleftrightarrow $ $f_{(n)}(x)$ exists}.
\]
\end{thm}

Marcinkiewicz and Zygmund proved their theorem by defining the sequence of differences $\widetilde{\Delta}_1(x,h)=f(x+h)-f(x)$ and $\widetilde{\Delta}_k(x,h)=\widetilde{\Delta}_{k-1}(x,2h)-2^{k-1}\widetilde{\Delta}_{k-1}(x,h)$ for $k\geq 2$, and showing that $\widetilde{\Delta}_n(x,h)$ is a scalar multiple of an $n$-th generalized Riemann difference whose associated generalized derivative satisfies their theorem. They were not interested in replacing 2 by other numbers, nor in other variants of their derivative, as the present paper considers, nor do they give an expression for the coefficients of their generalized Riemann derivative.

The goal of this paper is to determine a most inclusive generalization of Theorem MZ that addresses the above three questions that the classical theorem did not answer.

\subsection*{Results.}
We introduce three new kinds of $n$-th generalized Riemann derivatives, all depending on a real parameter $q$, with $q\neq 0,\pm 1$. The first two kinds are $q$-analogues of the $n$-th Riemann derivative $D_nf(x)$ and we call them $n$-th \emph{Gaussian Riemann derivatives}. These are the $n$-th generalized Riemann derivatives 
\[
\begin{aligned}
&\text{${_q}{D}_nf(x)$, based at $x,x+h,x+qh,x+q^2h,\ldots ,x+q^{n-1}h$},\\
&\text{and ${_q}{\bar D}_nf(x)$, based at $x+h,x+qh,x+q^2h,\ldots ,x+q^{n}h$.}
\end{aligned}
\]
The third, or the $n$-th \emph{symmetric Gaussian Riemann derivative}, is a $q$-analogue of the $n$-th symmetric Riemann derivative ${D}_n^sf(x)$. This is the $n$-th generalized Riemann derivative
\[
\text{${_q}{D}_n^sf(x)$, based at $(x),x\pm h,x\pm qh,x\pm q^2h,\ldots ,x\pm q^{m-1}h$,}
\]
where $m=\lfloor (n+1)/2\rfloor $ and~$(x)$ means that $x$ is taken only for $n$ even.
Note that $\widetilde{D}_n={_2}{D}_n$, so that ${_q}{D}_n$ is a generalization of $\widetilde{D}_n$, and ${_q}{\bar D}_n$ and ${_q}{D}_n^s$ are variants of this generalization.
The terms Gaussian and symmetric Gaussian will also refer to the differences, ${_q}{\Delta}_n(x,h;f)$, ${_q}\bar{\Delta}_n(x,h;f)$, and ${_q}{\Delta}_n^s(x,h;f)$, associated to these derivatives.

In section~1 we provide the explicit expressions for these differences in terms of the Gaussian binomial coefficients and determine the recursive relations for each of the three classes of differences, a generalization of the recursive relations found by~\cite{MZ}.

\medskip
Our main result is the following generalization of Theorem~MZ, proved in section~\ref{S2}, followed by a conjecture predicting that the result of the theorem is the best possible:

\begin{thmA}
Let $q$ be a real number with $q\neq 0,\pm 1$, and let $n$ be a positive integer. Then, for each function $f$ and point $x$,
\begin{enumerate}
\item[(i)\,] 
\text{both $f_{(n-1)}(x)$ and one of ${_q}D_nf(x)$ or ${_q}\bar D_nf(x)$ exist $\; \Longleftrightarrow \; $ $f_{(n)}(x)$ exists.}
\item[(ii)] 
\text{both $f_{(n-2)}^s(x)$ and ${_q}D_n^sf(x)$ exists $\; \Longleftrightarrow \; $ $f_{(n)}^s(x)$ exists.}
\end{enumerate}
\end{thmA}
When $n=1$, part (ii) of Theorem~A is still valid by ignoring the term $f_{(n-2)}^s(x)$ that does not make sense.
\begin{conjA}
{\rm
Fix a positive integer $n$ and let $D_{\mathcal{A}}$ be an $n$-th generalized Riemann differentiation without excess. Then:
\begin{enumerate}
\item[(i)\,]  If $D_{\mathcal{A}}$ is not Gaussian, then~for~all~$f$~and~$x$,
\[\qquad\qquad\text{both $f_{(n-1)}(x)$ and $D_{\mathcal{A}}f(x)$ exist $\;\;\; \not\hspace{-.05in}\Longrightarrow \; $ $f_{(n)}(x)$ exists.}\]

\item[(ii)] If $n\geq 3$ and $D_{\mathcal{A}}$ is symmetric but not Gaussian symmetric, then for all~$f$ and~$x$,
\[\qquad\qquad\text{both $f_{(n-2)}^s(x)$ and $D_{\mathcal{A}}f(x)$ exists $\;\;\; \not\hspace{-.05in}\Longrightarrow \; $ $f_{(n)}^s(x)$ exists.}\]
\end{enumerate}
}
\end{conjA}

Part~(i) of Conjecture~A is shown to be true in Proposition~\ref{P2.4}~for~$n=1$, and in Proposition~\ref{P2.5} for $n=2$, leaving it open for $n\geq 3$.
Part (ii) is easily proved false for $n=1$ and 2, and is proved true for $n=3$ and~4 in Proposition~\ref{P2.6}, leaving~it~open~for~$n\geq 5$.

\medskip
Conjecture~A is more general than the following conjecture, addressed in section~3, on Peano and Riemann derivatives, two of the most known generalized derivatives:

\begin{conjB}[\cite{AC2}]
{\rm For all functions $f$ and points $x$,
\begin{enumerate}
\item[(i)\,] When $n\geq 3$, $f_{(n-1)}(x)$ and $D_nf(x)$ exist $\;\;\not\hspace{-.04in}\Longrightarrow $ $f_{(n)}(x)$ exists.
\item[(ii)] When $n\geq 5$, $f_{(n-2)}^s(x)$ and $D_n^sf(x)$ exist $\;\;\not\hspace{-.04in}\Longrightarrow $ $f_{(n)}^s(x)$ exists.
\end{enumerate}
}
\end{conjB}

Both parts of this conjecture were first formally stated for $n\geq 3$ in \cite{AC2} as Conjectures 4.2 and 4.1.
Part (i) of Conjecture~B was proved for $n=3$ in \cite[Theorem~1(ii)]{ACF} via a clever example that does not extend to higher $n$; we prove Conjecture~B(i) for $n=7$ in Theorem~\ref{T3.2A}.
We update part (ii) of Conjecture~B to $n\geq 5$, based on us proving the asserted result false for~$n=3$ and 4 in Theorem~\ref{T3.2}(i). In addition, Theorem~\ref{T3.2}(ii) proves Conjecture~B(ii) for $n=5,6,7,8$, leaving it open for $n\geq 9$.

\subsection*{Details}
Let $q$ be a real number, with $q\neq 0,\,\pm 1$, and let $n$ be a positive integer. The quantum integer $n$, the quantum $n$ factorial, and the quantum $n$-choose-$k$ are, respectively,
\[
\text{$[n]_q=1+q+\cdots +q^{n-1}$, $[n]_q!=[1]_q[2]_q\cdots [n]_q$, and {\small $\left[ \begin{array}{c}n\\k\end{array}\right]_q=\frac {[n]_q!}{[n-k]_q!\cdot [k]_q!}$},}
\]
for $k=0,1,\ldots ,n$, where $[0]_q!=1$ and {\tiny $\left[ \begin{array}{c}n\\k\end{array}\right]_q$}$=0$ for $k>n$. Taking the limit as $q\rightarrow 1$ in $[n]_q$, $[n]_q!$, and {\tiny $\left[ \begin{array}{c}n\\k\end{array}\right]_q$}, respectively leads to $n$, $n!$, and $n\choose k$.
The Gaussian or $q$-binomial coefficients {\tiny $\left[ \begin{array}{c}n\\k\end{array}\right]_q$} are related to the Gaussian or $q$-binomial formula,
\begin{equation}\label{eq1}{\small
\text{$(a-b)(a-bq)(a-bq^2)\ldots (a-bq^{n-1})=\sum_{k=0}^n(-1)^kq^{{k}\choose 2}\left[ \begin{array}{c}n\\k\end{array}\right]_qa^{n-k}b^k,$}}
\end{equation}
whose limit as $q\rightarrow 1$ is the classical (Newton's) binomial formula.

To get an idea of how the quantum integers and Gaussian binomial coefficients are involved in the expressions of the $n$-th Gaussian Riemann differences,
\[
\begin{aligned}
&\text{${_q}{\Delta}_n(x,h;f)$, based at $x,x+h,x+qh,a+q^2h,\ldots ,x+q^{n-1}h$,}\\
&\text{${_q}\bar{\Delta}_n(x,h;f)$, based at $x+h,x+qh,a+q^2h,\ldots ,x+q^{n}h$, and}\\
&\text{${_q}{\Delta}_n^s(x,h;f)$, based at $(x),x\pm h,x\pm qh,x\pm q^2h,\ldots ,x\pm q^{m-1}h$, $m=\lfloor \tfrac {n+1}2\rfloor $,}
\end{aligned}
\]
using the Vandermonde relations, we computed the first two differences for $n=1,2,3$, and the third for $n=1,\ldots ,6$. After factoring their dominant coefficients, $\lambda_n$, $\bar\lambda_n$ and~$\lambda_n^s$, we observed that the remaining factors, ${_q}\widetilde{\Delta}_n(x,h;f)$, ${_q}\widetilde{\bar\Delta}_n(x,h;f)$ and ${_q}\widetilde{\Delta}_n^s(x,h;f)$, are differences with coefficients polynomials in $q$ related to $q$-binomial coefficients. 

Another observation we made after the computation of these differences for small $n$ was that the formulas for ${_q}{\Delta}_n(x,h;f)$, $n=1,2,3$, resemble the formulas for ${_q}{\Delta}_n^s(x,h;f)$, for $n=1,3,5$, and the formulas for~${_q}\bar{\Delta}_n(x,h;f)$, $n=1,2,3$, resemble the formulas for~${_q}{\Delta}_n^s(x,h;f)$, $n=2,4,6$. In this way,~it makes sense to have two $q$-analogues of the $n$-th forward Riemann difference for each~$n$,~so that the symmetric Gaussian case is a symmetric analogue of the (forward) Gaussian case.

The expressions for ${_q}{\Delta}_n(x,h;f)$, ${_q}{\bar\Delta}_n(x,h;f)$, and ${_q}{\Delta}_n^s(x,h;f)$, are given in Lemmas~\ref{L1}--\ref{L4}.
Recursive formulas for these differences are given by (\ref{eq:R1}), (\ref{eq7B}), and (\ref{eq11}), and the recursive relations for their difference quotients are given by (\ref{eq:R2}), (\ref{eq8B}), and (\ref{eq12}). The last three relations are not being used anywhere else in the paper; we included them for those readers interested to learn more about the combinatorics of the Gaussian differences.

The recursive relations obtained by Marcinkiewicz and Zygmund are crucial in the proof of Theorem~MZ, and, in the same way, our recursive relations, (\ref{eq:R1}), (\ref{eq7B}), and (\ref{eq11}), will be crucial in the proof of Theorem~A.
However, there is a difference in the way we obtained our differences, as we describe next:

The recursive relations in \cite{MZ} are given by definition, that is, by guessing, which most likely came after computing a number of differences for small $n$. From the recursive relations for the differences, one can deduce the recursive relations of their coefficients, and then realize that these are specializations of the Pascal's triangle identities for $q$-binomial coefficients, which in turn would have given the exact expressions for the MZ-differences.
Our equivalent method of obtaining these differences goes the other way around. Since the independent computational goal of the paper was the deduction of the expressions for the Gaussian differences, we used the intuition gained from the computation of differences for small $n$ towards guessing the exact specializations of the q-binomial formula that are equivalent to the expressions of the three differences. Then, the resulting Pascal's triangle identities for the known coefficients led to the recursive relations for the same differences.
\[
\ast\quad\ast\quad\ast
\]

The equivalence between generalized derivatives is an almost a century old problem. It was initiated by Kintchine in \cite{Ki} (1927), who proved that the first symmetric derivative is a.e. equivalent to the first Peano derivative. This was greatly extended by Marcinkiewicz and Zygmund in \cite{MZ} (1936) to the a.e. equivalence between the $n$-th Peano and the $n$-th symmetric Riemann derivatives, and further by Ash in \cite{As} (1967) who showed that any $n$-th generalized Riemann derivative of a function $f$ is a.e. equivalent to the $n$-th Peano derivative on a measurable set.
In particular, any two generalized Riemann derivatives of~$f$ of the same order are a.e. equivalent on a measurable set. The pointwise equivalence between Peano and generalized Riemann differentiation is studied in \cite{ACCs}; an application of this to continuity is found in \cite{AAC}. The pointwise equivalence between generalized Riemann and ordinary differentiations for continuous functions is addressed in~\cite{C}. Pointwise equivalences and pointwise implications between any two generalized Riemann derivatives of a real or complex function $f$ are investigated in \cite{ACCh} and \cite{ACCH}. In particular, the above mentioned single equivalent class breaks up into numerous smaller equivalence classes, and these are described explicitly. Quantum Riemann derivatives were introduced in \cite{AC,ACR}. Multidimensional Riemann derivatives are explored in \cite{AC1}. These recent articles have shown numerous connections between generalized Riemann derivatives and linear and abstract algebra, recursive set theory, symmetric functions, complex and numerical analysis. For more on Peano and generalized Riemann differentiation, see \cite{As1,ACF1,CF,F,F1,FR,GR,LPW,RAA2}.

Generalized Riemann derivatives have many applications in the theory of trigonometric series \cite{SZ,Z}. They were shown to satisfy properties similar to those for ordinary derivatives, such as convexity, monotonicity, and the mean value theorem \cite{AJ,FFR,GGR,HL1,HL2,MM,T1,W}. Surveys on generalized derivatives are found in \cite{As2} and \cite{EW}.
\section{Explicit formulas for Gaussian Riemann differences}\label{S1}
In this section we provide the explicit formulas for the $n$-th Gaussian differences ${_q}{\Delta}_n(x,h;f)$, ${_q}{\bar\Delta}_n(x,h;f)$, and ${_q}{\Delta}_n^s(x,h;f)$. These are proved by reference to various specializations of the Gaussian binomial formula.

\subsection{Forward Gaussian Riemann differences.} 
Taking~$n-1$ instead of $n$ in (\ref{eq1}) and setting $b=q$ leads to the Gaussian binomial formula
\begin{equation}\label{eq2}\small
(a-q)(a-q^2)\ldots (a-q^{n-1})=\sum_{k=0}^{n-1}(-1)^kq^{{k+1}\choose 2}\left[ \begin{array}{c}n-1\\k\end{array}\right]_qa^{n-1-k}.
\end{equation}
When $a=1$, equation (\ref{eq2}) is equivalent to
\begin{equation}\label{eq3}\small
\sum_{k=0}^{n-1}(-1)^kq^{{k+1}\choose 2}\left[ \begin{array}{c}n-1\\k\end{array}\right]_q-(1-q)(1-q^2)\ldots (1-q^{n-1})=0;
\end{equation}
when $a=q^j$, for $j=1,\ldots ,n-1$, the same equation is
\begin{equation}\label{eq4}\small
\sum_{k=0}^{n-1}(-1)^kq^{{k+1}\choose 2}\left[ \begin{array}{c}n-1\\k\end{array}\right]_q\left(q^{n-1-k}\right)^j=0;
\end{equation}
and when $a=q^n$, the same equation becomes
\begin{equation}\label{eq5}\small
\sum_{k=0}^{n-1}(-1)^kq^{{k+1}\choose 2}\left[ \begin{array}{c}n-1\\k\end{array}\right]_q\left(q^{n-1-k}\right)^n=(q^n-q)(q^n-q^2)\ldots (q^n-q^{n-1}).
\end{equation}

The following lemma provides the expression of the $n$-th Gaussian Riemann difference.

\begin{lemma}\label{L1}
The $n$-th forward Gaussian Riemann difference has the expression
{\small \[
{_q}{\Delta}_n(x,h;f)=\lambda_n\cdot {_q}\widetilde{\Delta}_n(x,h;f),
\]}
where
{\small \[
{_q}\widetilde{\Delta}_n(x,h;f)=\sum_{k=0}^{n-1}(-1)^kq^{{k+1}\choose 2}\left[ \begin{array}{c}n-1\\k\end{array}\right]_qf(x+q^{n-1-k}h)-(1-q)(1-q^2)\ldots (1-q^{n-1})f(x)
\]}
and {\small $\lambda_n=n!/\left((q^n-q)(q^n-q^2)\ldots (q^n-q^{n-1})\right)$.}
\end{lemma}

\begin{proof}
By (\ref{eq3}), (\ref{eq4}), (\ref{eq5}) the difference {\small $\lambda_n\cdot {_q}\widetilde{\Delta}_n(x,h;f)$} based at $x,x+h,x+qh,\ldots ,x+q^{n-1}h$ satisfies the $n$-th Vandermonde relations, so it must equal ${_q}{\Delta}_n(x,h;f)$.
\end{proof}

Using the quantum Pascal triangle identity
{\small\begin{equation}\label{Pascal}
\left[ \begin{array}{c}n-1\\k\end{array}\right]_q=\left[ \begin{array}{c}n-2\\k\end{array}\right]_q+q^{n-1-k}\left[ \begin{array}{c}n-2\\k-1\end{array}\right]_q,
\end{equation}}and the expressions for the differences {\small $\widetilde{\Delta}_n(x,h;f)$} given in Lemma~\ref{L1}, one can easily show that the same differences can also be computed recursively as
\begin{equation}\label{eq:R1}\small
\begin{aligned}
{_q}\widetilde{\Delta}_{1}\left(  x,h;f\right)    & =f\left(  x+h\right)-f\left(  x\right) , \\
{_q}\widetilde{\Delta}_{n}\left(  x,h;f\right)    & ={_q}\widetilde{\Delta}_{n-1}\left(  x,qh;f\right)  -q^{n-1}\cdot {_q}\widetilde{\Delta}_{n-1}\left(
x,h;f\right)  ,\qquad (n\geq 2).
\end{aligned}
\end{equation}

Another recursive way to define the $n$-th Gaussian Riemann difference is by using the sequence of generalized Riemann difference quotients given by
{\small \[
{_q}\widetilde{D}_n(x,h;f)={_q}{\Delta}_{n}\left(  x,h;f\right) /h^n.
\]}
By (\ref{eq:R1}) and the expression of $\lambda_n$ in Lemma~\ref{L1}, one can prove by induction on $n$ that this sequence satisfies the following recursive relation:

\begin{equation}\label{eq:R2}\small
\begin{aligned}
{_q}\widetilde{D}_1(x,h;f)&=\frac {f(x+h)-f(x)}h,\\
{_q}\widetilde{D}_n(x,h;f)&=n\frac {{_q}\widetilde{D}_{n-1}(x,qh;f)-{_q}\widetilde{D}_{n-1}(x,h;f)}{(q^{n-1}-1)h} ,\qquad (n\geq 2).
\end{aligned}
\end{equation}
As an $n$-th generalized Riemann difference quotient, the same sequence enjoys the property that for each $n$, {\small ${_q}\widetilde{D}_n(x,h;f)=f^{(n)}(x)$}, for all polynomials $f$ of degree~$\leq n$.

The \emph{$n$-th Gaussian Riemann derivative} of a function $f$ at $x$ is the $n$-th generalized Riemann derivative
{\small \[
{_q}D_nf(x):=\lim_{h\rightarrow 0}{_q}\widetilde{D}_n(x,h;f)=\lim_{h\rightarrow 0}\frac {{_q}\Delta_n(x,h;f)}{h^n}=\lim_{h\rightarrow 0}\frac {\lambda_n\cdot {_q}\widetilde{\Delta}_n(x,h;f)}{h^n}.
\]}It is a $q$-analogue of the $n$-th (forward) Riemann derivative, but not the only one. For example, it is different from the $n$-th quantum Riemann derivative defined in~\cite{ACR},
which satisfies $q$-Vandermonde relations instead of ordinary Vandermonde relations, hence is not an $n$-th generalized Riemann derivative.

Another $q$-analogue of the $n$-th Riemann derivative which \emph{is} an $n$-th generalized Riemann derivative is the unique $n$-th generalized Riemann derivative based at $x+h,x+qh,\ldots ,x+q^nh$ whose expression is given explicitly in the following~lemma:

\begin{lemma}\label{L2}
The $n$-th generalized Riemann difference based at $x+h,x+qh,\ldots ,x+q^nh$ has the expression
{\small \[
{_q}{\bar\Delta}_n(x,h;f)=\bar\lambda_n\sum_{k=0}^n(-1)^kq^{k\choose 2}\left[ \begin{array}{c}n\\k\end{array}\right]_qf(x+q^{n-k}h),
\]}
where {\small $\bar\lambda_n=n!/\left((q^n-1)(q^n-q)\ldots (q^n-q^{n-1})\right)$}.
\end{lemma}

\begin{proof}
As in Lemma~\ref{L1}, the expression is implied by the $q$-binomial formula
{\small \[
(a-1)(a-q)(a-q^2)\ldots (a-q^{n-1})=\sum_{k=0}^n(-1)^kq^{{k}\choose 2}\left[ \begin{array}{c}n\\k\end{array}\right]_qa^{n-k}
\]}obtained from (\ref{eq1}) by taking $b=1$,
since its Vandermonde relations are deduced from this formula, by taking $a=q^j$ for $j=0,1,\ldots ,n$.
\end{proof}

The identity ({\ref{Pascal}) for $n$ instead of $n-1$ and Lemma~\ref{L2} can be employed to deduce the recursive relations between the Gaussian differences {\small ${_q}{\bar\Delta}_n(x,h;f)$}. These are
\begin{equation}\label{eq7B}\small
\begin{aligned}
{_q}\widetilde{\bar\Delta}_{1}\left(  x,h;f\right)    & =f\left(  x+qh\right)-f\left(  x+h\right) , \\
{_q}\widetilde{\bar\Delta}_{n}\left(  x,h;f\right)    & ={_q}\widetilde{\bar\Delta}_{n-1}\left(  x,qh;f\right)  -q^{n-1}\cdot {_q}\widetilde{\bar\Delta}_{n-1}\left(
x,h;f\right)  ,\qquad (n\geq 2),
\end{aligned}
\end{equation}
which in turn can be involved in deducing the following recursive relations satisfied by the defining difference quotients
{\small ${_q}\widetilde{\bar D}_n(x,h;f):={_q}{\bar\Delta}_{n}\left(  x,h;f\right)/h^n$}:

\begin{equation}\label{eq8B}\small
\begin{aligned}
{_q}\widetilde{\bar D}_1(x,h;f)&=\frac {f(x+qh)-f(x+h)}{(q-1)h},\\
{_q}\widetilde{\bar D}_n(x,h;f)&=n\frac {{_q}\widetilde{\bar D}_{n-1}(x,qh;f)-{_q}\widetilde{\bar D}_{n-1}(x,h;f)}{(q^{n}-1)h},\qquad (n\geq 2).
\end{aligned}
\end{equation}

The $n$-th (forward) Gaussian Riemann derivative ${_q}{\bar D}_nf(x)$ is defined by the limit:
{\small \[
{_q}{\bar D}_nf(x)=\lim_{h\rightarrow 0}{_q}\widetilde{\bar D}_n(x,h;f)=\lim_{h\rightarrow 0}{_q}{\bar \Delta}_n(x,h;f)/h^n.
\]}

\subsection{Symmetric Gaussian Riemann differences.} Throughout this section, $n$ will be a fixed positive integer and $m=\lfloor (n+1)/2\rfloor $. Recall from the introduction that the expressions of the symmetric Gaussian Riemann differences depend on the parity of $n$. Their exact formulas for general $n$ will be deduced from two specializations of the $q$-binomial formula in a similar way as we did for the forward Gaussian Riemann derivatives.
The first of these specializations,
\begin{equation}\label{eq8}\small
(a-q^2)(a-q^4)\ldots (a-q^{2(m-1)})=\sum_{k=0}^{m-1}(-1)^kq^{k(k+1)}\left[ \begin{array}{c}m-1\\k\end{array}\right]_{q^2}a^{m-1-k},
\end{equation}
is obtained from (\ref{eq2}) by replacing $n$ with $m$ and $q$ with $q^2$.

The following lemma provides the expression for the $n$-th symmetric Gaussian Riemann derivative in the case when $n$ is even.

\begin{lemma}\label{L3}
When $n=2m$, the $n$-th symmetric Gaussian Riemann difference is the $n$-th generalized Riemann difference based at {\small $x,x\pm h,x\pm qh,x\pm q^2h,\ldots ,x\pm q^{m-1}h$}. This has the expression {\small ${_q}{\Delta}_n^s(x,h;f)=\lambda_n^s\cdot {_q}\widetilde{\Delta}_n^s(x,h;f)$}, where
{\small \[\begin{aligned}
{_q}\widetilde{\Delta}_n^s(x,h;f)&=\sum_{k=0}^{m-1}(-1)^kq^{k(k+1)}\left[ \begin{array}{c}m-1\\k\end{array}\right]_{q^2}\left\{f(x+q^{m-1-k}h)+f(x-q^{m-1-k}h)\right\}\\
&-2(1-q^2)(1-q^4)\ldots (1-q^{2(m-1)})f(x)
\end{aligned}
\]}and {\small $\lambda_n^s=n!/\left(2(q^n-q^2)(q^n-q^4)\ldots (q^n-q^{n-2})\right)$}.
\end{lemma}

\begin{proof}
The Vandermonde conditions for $j$ even are verified in the same way as in the proof of Lemma~\ref{L1}, this time using (\ref{eq3}), (\ref{eq4}) and (\ref{eq5}) with $n$ replaced by $m$ and~$q$ replaced by $q^2$. For odd $j$ and even difference, the Vandermonde conditions are trivially satisfied.
\end{proof}

The second specialization of the $q$-binomial formula that will be needed in dealing with symmetric Gaussian Riemann differences is
\begin{equation}\label{eq9}\small
(a-q)(a-q^3)\ldots (a-q^{2m-3})=\sum_{k=0}^{m-1}(-1)^kq^{k^2}\left[ \begin{array}{c}m-1\\k\end{array}\right]_{q^2}a^{m-1-k}.
\end{equation}

The following lemma provides the expression for the $n$-th symmetric Gaussian Riemann derivative in the case when $n$ is odd.

\begin{lemma}\label{L4}
When $n=2m+1$, the $n$-th symmetric Gaussian Riemann difference is the $n$-th generalized Riemann difference based at {\small $x\pm h,x\pm qh,x\pm q^2h,\ldots ,x\pm q^{m-1}h$}. This has the expression {\small ${_q}{\Delta}_n^s(x,h;f)=\lambda_n^s\cdot {_q}\widetilde{\Delta}_n^s(x,h;f)$}, where
{\small \[
\widetilde{\Delta}_n^s(x,h;f)=\sum_{k=0}^{m-1}(-1)^kq^{k^2}\left[ \begin{array}{c}m-1\\k\end{array}\right]_{q^2}\left\{f(x+q^{m-1-k}h)-f(x-q^{m-1-k}h)\right\}
\]}and {\small $\lambda_n^s=n!/\left(2(q^n-q)(q^n-q^3)\ldots (q^n-q^{n-2})\right)$}.
\end{lemma}

\begin{proof}
The Vandermonde conditions for $j$ odd are verified in the same way as in Lemma~\ref{L1}, this time using (\ref{eq3}), (\ref{eq4}), (\ref{eq5}) with $n$ replaced by $m$ and $q$ replaced by~$q^2$. The Vandermonde conditions are trivially satisfied for even~$j$ and odd difference.
\end{proof}

Since both expressions in Lemmas~\ref{L3} and \ref{L4} involve the same $q$-binomial coefficients, the Pascal triangle identity (\ref{Pascal}) with $m$ instead of $n$ and $q^2$ instead of $q$ can be used to inductively deduce the following combined recursive relation for all symmetric Gaussian Riemann differences:
\begin{equation}\label{eq11}\small
\begin{aligned}
{_q}\widetilde{\Delta}_{1}^s\left(  x,h;f\right)    & =f\left(  x+h\right)-f\left(  x-h\right),  \\
{_q}\widetilde{\Delta}_{2}^s\left(  x,h;f\right)    & =f\left(  x+h\right)-2f(x)+f\left(  x-h\right),  \\
{_q}\widetilde{\Delta}_{n}^s\left(  x,h;f\right)    & ={_q}\widetilde{\Delta}_{n-1}^s\left(  x,qh;f\right)  -q^{n-2}\cdot {_q}\widetilde{\Delta}_{n-2}^s\left(
x,h;f\right)  ,\qquad (n\geq 3).
\end{aligned}
\end{equation}

Finally, we can use the recursive relations (\ref{eq11}) and the expressions for $\lambda_n^s$ provided by Lemmas~\ref{L3} and~\ref{L4} to inductively prove that the $n$-th symmetric Gaussian Riemann quotients
{\small \[
{_q}\widetilde{D}_n^s(x,h;f):={_q}{\Delta}_{n}^s\left(  x,h;f\right) /h^n=\lambda_n^s\cdot {_q}\widetilde{\Delta}_n^s(x,h;f)/h^n
\]}
satisfy the following recursive relations:
\begin{equation}\label{eq12}\small
\begin{aligned}
{_q}\widetilde{D}_1^s(x,h;f)&=\frac {f(x+h)-f(x)}h,\;
{_q}\widetilde{D}_2^s(x,h;f)=\frac {f\left(  x+h\right)-2f(x)+f\left(  x-h\right)}{h^2},\\
{_q}\widetilde{D}_n^s(x,h;f)&=n(n-1)\frac {{_q}\widetilde{D}_{n-2}^s(x,qh;f)-{_q}\widetilde{D}_{n-2}^s(x,h;f)}{(q^{n-2+(n\mod 2)}-1){h^2}},\qquad (n\geq 3).
\end{aligned}
\end{equation}
The \emph{$n$-th symmetric Gaussian Riemann derivative} of a function $f$ at $x$ is the $n$-th generalized Riemann derivative
{\small \[
{_q}D_n^sf(x):=\lim_{h\rightarrow 0}{_q}\widetilde{D}_n^s(x,h;f).
\]}
It is a $q$-analogue of the $n$-th symmetric Riemann derivative, and is different from the $n$-th quantum symmetric Riemann derivative defined in~\cite{AC}.

\section{Proof of Theorem A and Evidence for Conjecture A}\label{S2}
This section has two parts: the proof of Theorem~A; and the evidence for Conjecture~A.

\subsection{Proof of Theorem~A} We are now ready to proceed with the proof of Theorem~A. For this, we will need the following lemma:

\begin{lemma}\label{L2.1}
Let $q$ be a real number with $q\neq 0,\pm 1$, and let $n$ be a positive integer. Then for each function $f$ and point $x$,
\begin{enumerate}
\item[(i)\,] 
\text{${_q}D_nf(x)$ exists $\; \Longleftrightarrow \; $ ${_{q^{-1}}}D_nf(x)$ exists.}
\item[(ii)] 
\text{${_q}D_n^sf(x)$ exists $\; \Longleftrightarrow \; $ ${_{q^{-1}}}D_n^sf(x)$ exists.}
\end{enumerate}
\end{lemma}

\begin{proof} Part (i) follows from the Gaussian Riemann differences ${_{q^{-1}}}\Delta_n(x,h;f)$ and ${_{q}}\Delta_n(x,h;f)$ being scales of each other by $q^{\pm (n-1)}$, since they are respectively based at $x,x+h,x+q^{-1}h,\ldots ,x+q^{-n+1}h$ and~$x,x+h,x+qh,\ldots ,x+q^{n-1}h$. Part~(ii) follows from the similar property between ${_{q^{-1}}}\Delta_n^s(x,h;f)$ and ${_{q}}\Delta_n^s(x,h;f)$.
\end{proof}

For simplicity, in the proof of Theorem~A we denote a difference $\Delta (0,h;f)$~as~$\Delta (h)$.

\begin{proof}[Proof of Theorem~A]
(i) As the existence of the $n$-th Peano derivative {\small $f_{(n)}(x)$} both assumes the existence of every lower order Peano derivatives of $f$ at $x$ and implies every $n$-th generalized Riemann derivative of $f$ at $x$, the reverse implication is clear.

Conversely, suppose that both $f_{(n-1)}(x)$ and ${_q}D_nf(x)$ exist. By Lemma~\ref{L2.1}, we may assume that $|q|>1$. And eventually by translating the graph of $f$ to the left by $x$ we may assume that $x=0$, and by subtracting from $f$ a degree $n$ polynomial we may further assume that {\small $f_{(n-1)}(0)=0$} and~{\small ${_q}D_nf(0)=0$}, or {\small ${_q}\Delta_n(h)=o(h^n)$}. This is equivalent to having~{\small ${_q}\widetilde{\Delta}_n(h)=o(h^n)$} since {\small ${_q}\lambda_n$} is independent of $h$. The last equality means that
for each $\varepsilon >0$, there is a $\delta >0$ such that {\small $|h|<\delta \Rightarrow |{_q}\widetilde{\Delta}_n(h)|<\varepsilon |h|^n$}. Then, by (\ref{eq:R1}), we deduce
{\small\[
\left|{_q}\widetilde{\Delta}_{n-1}(qh)-q^{n-1}{_q}\widetilde{\Delta}_{n-1}(h)\right|< \varepsilon |h|^n,\quad\left|{_q}\widetilde{\Delta}_{n-1}(h)-q^{n-1}{_q}\widetilde{\Delta}_{n-1}\left(\frac hq\right)\right|< \varepsilon \left|\frac hq\right|^n,\ldots 
\]}{\small\[
\ldots ,\left|{_q}\widetilde{\Delta}_{n-1}\left(\frac h{q^{k-1}}\right)-q^{n-1} {_q}\widetilde{\Delta}_{n-1}\left(\frac h{q^k}\right)\right|< \varepsilon \left|\frac h{q^k}\right|^n.
\]}We multiply these inequalities resp. by {\small $1, q^{n-1},q^{2(n-1)},\ldots ,q^{k(n-1)}$} and~add.~The triangle inequality makes the left side telescope, while the right side is a geometric series. Then
{\small\[
\left|{_q}\widetilde{\Delta}_{n-1}(qh)-q^{(k+1)(n-1)}\cdot {_q}\widetilde{\Delta}_{n-1}\left(\frac h{q^k}\right)\right|< \frac {q}{q-1}\cdot \varepsilon \left| h\right|^n.
\]}
The second term on the left can be neglected, since it is
{\small
$(qh)^{n-1}{\left({_q}\widetilde{\Delta}_{n-1}(h/{q^k})\right)}/{\left(h/{q^k}\right)^{n-1}}$}
and this approaches $0$ as $k\rightarrow \infty $, by the hypothesis
{\small $f_{(n-1)}(0)={_q}D_nf(0)=0$}.~Therefore,
{\small\[
\left|{_q}\widetilde{\Delta}_{n-1}(qh)\right|< \frac q{q-1}\cdot \varepsilon |h|^n\text{, that is, }\left|{_q}\widetilde{\Delta}_{n-1}(h)\right|=o(h^n).
\]}
By the independence on $h$ of ${_q}\lambda_{n-1}$, this is equivalent to $\left|{_q}{\Delta}_{n-1}(h)\right|=o(h^n)$.
Similarly, one can deduce that $\left|{_q}\Delta_{n-2}(h)\right|=o(h^n)$, and so on. At the end, $\left|{_q}\Delta_{1}(h)\right|=o(h^n)$ means that $f(h)-f(0)=f(h)=o(h^n)$, and hence $f_{(n)}(0)=0$, as needed.
The direct implication under the hypothesis that both $f_{(n-1)}(x)$ and ${_q}\bar D_nf(x)$ exist  is proved along the same lines. The proof of part (ii) is similar to the proof of part~(i).
\end{proof}

\subsection{Evidence for Conjecture~A}
The remaining part of the section analyzes the evidence towards this conjecture by proving its asserted result in a number of~cases.

\medskip
The following proposition shows that part (i) of Conjecture~A is true for $n=1$.

\begin{proposition}\label{P2.4}
Let $D_{\mathcal{A}}f(x)$ be a first generalized Riemann derivative without excess which is not a Gaussian Riemann derivative. Then:
\begin{enumerate}
\item[(i)\,] $D_{\mathcal{A}}f(x)=f_{(1)}^s(x)$;
\item[(ii)] both $f_{(0)}(x)$ and $D_{\mathcal{A}}f(x)$ exist $\;\;\; \not\hspace{-.05in}\Longrightarrow \; $ $f_{(1)}(x)$ exists, for all $f$ and $x$.
\end{enumerate}
\end{proposition}

\begin{proof}
(i) The hypothesis that $D_{\mathcal{A}}f(x)$ is a first generalized Riemann derivative without excess makes its difference $\Delta_{\mathcal{A}}f(x)=A_1f(x+a_1h)+A_2f(x+a_2h)$, for some $A_1,A_2,a_1,a_2$, with $a_1\neq a_2$. If $a_1a_2=0$, say $a_1=0$, then $\Delta_{\mathcal{A}}f(x)$ is the scale by~$a_2$ of the Riemann difference $\Delta_1f(x)=f(x+h)-f(x)={_q}\Delta_1f(x)$ for any~$q$, hence $D_{\mathcal{A}}f(x)$ is Gaussian, a contradiction.
If $a_1a_2\neq 0$ and $|a_1|\neq |a_2|$, then $\Delta_{\mathcal{A}}f(x)$ is the scale by~$a_1$ of the first difference based at $x+h,x+qh$ for $q=a_2/a_1$, which is ${\,}_q\bar\Delta_1f(x)$, hence $D_{\mathcal{A}}f(x)$ is Gaussian, a contradiction. In the remaining case $a_1=-a_2\neq 0$, $\Delta_{\mathcal{A}}(x,h;f)$ is the scale by $a_1$ of the symmetric difference $\Delta_1^s(x,h;f)$, hence $D_{\mathcal{A}}f(x)=f_{(1)}^s(x)$.

(ii) As an example, take the function $f(t)=|t-x|$, which is continuous at $t=x$ hence $f_{(0)}(x)=f(x)=0$, it has $f(x+h)-f(x-h)=0$ hence $D_{\mathcal{A}}f(x)=f_{(1)}^s(x)=0$, while $f_{(1)}(x)=\lim_{h\rightarrow 0} \{f(x+h)-f(x)\}/h=\lim_{h\rightarrow 0}|h|/h$ does not exist.
\end{proof}

The next proposition shows that part (i) of Conjecture~A is also true for $n=2$.

\begin{proposition}\label{P2.5}
Let $D_{\mathcal{A}}f(x)$ be a second generalized Riemann derivative without excess which is not a Gaussian Riemann derivative. Then:
\begin{enumerate}
\item[(i)\,] up to a scale, $\Delta_{\mathcal{A}}(x,h;f)$ is either based at $x\pm h,x+qh$, for $q\neq 0,\pm 1$, or based at $x+h,x+ph,x+qh$, where $p,q\neq 0,\pm 1$, $p\neq \pm q$ and none of $p$ and $q$ is the square of the other;
\item[(ii)] both $f_{(1)}(x)$ and $D_{\mathcal{A}}f(x)$ exist $\;\;\; \not\hspace{-.05in}\Longrightarrow \; $ $f_{(2)}(x)$ exists, for all $f$ and $x$.
\end{enumerate}
\end{proposition}

\begin{proof}
(i) The hypothesis that $D_{\mathcal{A}}f(x)$ is second order and without excess makes it have base points $x+a_1h,x+a_2h,x+a_3h$ for distinct $a_1,a_2,a_3$. If one of $a_1,a_2,a_3$ is zero, say $a_1=0$, then up to a scale by $a_2^{-1}$, $D_{\mathcal{A}}f(x)$ is based at $x,x+h,x+qh$, for $q=a_2/a_1$, that is $D_{\mathcal{A}}f(x)$ is Gaussian, a contradiction. Thus $a_1,a_2,a_3$ are all non-zero. If two of them add up to zero, say $a_1+a_2=0$, then up to a scale by $a_1^{-1}$, the difference has the first outlined form. Otherwise, a scale by $a_1^{-1}$ and discounting the Gaussian case makes the difference have the second outlined form.

(ii) Our example has $D_{\mathcal{A}}f(x)$ based at $x+h,x+2h,x+3h$, that is, $\Delta_{\mathcal{A}}(x,h;f)=f(x+3h)-2f(x+2h)+f(x+h)$. Eventually by shifting the graph of $f$ to the left by~$x$,~we may assume that $x=0$. Take $G=\langle 2,3\rangle =\{2^r3^s\mid r,s\text{ integers}\}$ and define $f$ as
\[
f(x)=(-1)^{r+s}x^2,\quad\text{ if }x=2^r3^s\in G,
\]
and $f(x)=0$, otherwise. Then $f(h)=o(h)$ as $h\rightarrow 0$, hence $f_{(0)}(0)=f_{(1)}(0)=0$, while $f_{(2)}(0)$ does not exist, due to $\lim_{h\rightarrow 0} f(h)/h^2=0,\pm 1$. Moreover, when $h=2^r3^s\in G$, $\Delta_{\mathcal{A}}(0,h;f)=f(3h)-2f(2h)+f(h)=(-1)^{r+s}(-3^2+2\cdot 2^2+1\cdot 1^2)h^2=0$, and when $h\notin G$, $\Delta_{\mathcal{A}}(0,h;f)=0-2\cdot 0+0=0$, and so $D_{\mathcal{A}}f(0)=0$.
\end{proof}

We now turn to the symmetric case addressed in part (ii) of Conjecture~A. The first symmetric Riemann derivative $D_1^sf(x)$ is up to a scale the only first symmetric generalized Riemann derivative of $f$ at $x$ without excess, and since its definition is the same as the definition of $f_{(1)}^s(x)$, the conjecture is false in the case $n=1$. Same story for $n=2$.

\medskip
The following proposition shows that Conjecture~A(ii) is true for $n=3$ and 4.

\begin{proposition}\label{P2.6}
Each order 3 or 4 symmetric generalized Riemann derivative without excess is a symmetric Gaussian Riemann derivative.
\end{proposition}

\begin{proof}
Let $D_{\mathcal{A}}f(x)$ be a symmetric generalized Riemann difference of order $n=3$ or 4. Then it is based at $(x),x\pm ph,x\pm qh$, for $0<p<q$, which then scaled by~$p^{-1}$ becomes a symmetric Gaussian Riemann difference.
\end{proof}

\section{Updating Conjecture B}

In this section we prove in Theorem~\ref{T3.2A} that Conjecture~B(i) is true for $n=7$, and we update Conjecture~B(ii) to $n\geq 5$ by disproving the asserted result for $n=3$ and~4 in Theorem~\ref{T3.2}(i), and positively answer the predicted result for $n=3,4,\ldots ,8$ in Theorem~\ref{T3.2}(ii), leaving it open for $n\geq 9$.

\medskip
The following theorem gives answers to Conjecture~B(i), for $n=7$.

\begin{theorem}\label{T3.2A}
When $n=7$, Conjecture~B(i) is true.
\end{theorem}

\begin{proof}
Let $G=\langle 2,3,5,7\rangle $ be the subgroup of the multiplicative group of non-zero real numbers generated by $2,3,5,7$, and let $f$ be the real function defined by the expression:
\[
f(h)=(-1)^{n+p}h^s,\text{ for }h=2^m3^n5^p7^q\in G,
\]
and $f(h)=0$, for $h\notin G$, where $s$ is a real number between 6 and 7 that will be determined later. Since $6<s<7$, $f$ is six times Peano differentiable at $x=0$, but not seven times. The difference $\Delta_7(h):=\Delta_7(0,h;f)$, for $h=2^m3^n5^p7^q\in G$, is
\[
\begin{aligned}
\Delta_7(h)&=f(7h)-7f(6h)+21f(5h)-35f(4h)+35f(3h)-21f(2h)+7f(h)-f(0)\\
&=(7^s+7\cdot 6^s-21\cdot 5^s-35\cdot 4^s-35\cdot 3^s-21\cdot 2^s+7)h^s,
\end{aligned}
\]
and $\Delta_7(h)=0$ for $h\notin G$.
Let $\varphi :[6,7]\rightarrow \mathbb{R}$ be the function defined by the expression in the above parenthesis, where the constant $s$ is replaced with a variable $x$. Then $\varphi (6)=-54,096<0$ and $\varphi (7)=489,804>0$, so, by continuity, $\varphi (s)=0$, for some $s$ between 6 and 7. For this $s$, $\Delta_7(h)=0$ for all $h$, hence $D_7f(0)=0$, proving the result.
\end{proof}

\medskip
The following theorem gives answers to Conjecture~B(ii), for $n=3,4,\ldots ,8$.

\begin{theorem}\label{T3.2}
The following are answers to Conjecture~B(ii) for small values of $n$.
\begin{enumerate}
\item[(i)\,] When $n=3$ or $4$, the conjecture is false.
\item[(ii)] When $n=5,6,7,8$, the conjecture is true.
\end{enumerate}
\end{theorem}

\begin{proof}
(i) When $n=3$ or 4, the symmetric Riemann derivative $D_n^sf(x)$ is symmetric Gaussian, either by Proposition~\ref{P2.6} or directly by observing that
{\small\[
\begin{aligned}
2{\Delta}_3^s(x,h;f)&=f(x+3h)-3f(x+h)+3f(x-h)-f(x-3h)=2\cdot {_3}{\Delta}_3^s(x,h;f),\\
{\Delta}_4^s(x,h;f)&=f(x+2h)-4f(x+h)+6f(x)-4f(x-h)+f(x-2h)={_2}{\Delta}_4^s(x,h;f).
\end{aligned}
\]}The result then follows from Theorem~A(ii).

(ii) When $n=5$, let $G=\langle 3,5\rangle =\{3^m5^n\mid m,n\in\mathbb{Z}\}$ be the multiplicative subgroup
of the rationals generated by 3 and 5, and let $f:\mathbb{R}\rightarrow \mathbb{R}$ be defined as
\[
{\small \text{$f(x)=(-1)^{m+n}x^k$}}, \text{ for }{\small \text{$x=3^m5^n\in G$,}}
\]
and $f(x)=0$ for $x\notin G$, where $k$, $3<k<4$, is to be determined. Compute the~expression
\[
{\small \text{$\frac 12\{f(h)-f(-h)\}=\frac 12\cdot (-1)^{m+n}\cdot (\pm |h|^k)$}},\quad\text{ for }{\small h\in G},
\]
and {\small $\frac 12\{f(h)-f(-h)\}=0$} for {\small $h\notin G$}, to deduce {\small $f_{(3)}(0)=0$}, since {\small $\frac 12\{f(h)-f(-h)\}=o(h^3)$}, and {\small $f_{(5)}(0)$} does not exist since {\small $\lim_{h\rightarrow 0}\frac 12\{f(h)-f(-h)\}/h^5$} does not exist. A scale by 2 of the difference~{\small $\Delta_5^s(0,h;f)$} is the difference {\small$2^{-5}\cdot \Delta_5^s(0,2h;f)$}, where the difference
{\small $\Delta_5^s(0,2h;f)=f(5h)-5f(3h)+10f(h)-10f(-h)+5f(-3h)-6f(-5h)
=(-1)^{m+n+1}5^k|h|^k-5\cdot (-1)^{m+n+1}3^k|h|^k+10\cdot (-1)^{m+n}|h|^k 
=(-1)^{m+n+1}(5^k-5\cdot 3^k-10)|h|^k$}.
Denote {\small $\varphi (k)=5^k-5\cdot 3^k-10$} and observe that {\small $\varphi (3)\varphi (4)<0$}, and so by continuity, {\small $\varphi (k)=0$} for some $k$ in the open interval $(3,4)$. Then the $f$ for this $k$ has {\small $D_5^sf(0)=0$}.

When $n=6$, let $G=\langle 2,3\rangle =\{2^m3^n\mid m,n\in\mathbb{Z}\}$ and take $f$ to be the function
\[
{\small \text{$f(x)=(-1)^{m+n}x^k$}}, \quad\text{ for }{\small \text{$x=2^m3^n\in G$,}}
\]
and $f(x)=0$, for $x\notin G$, where the real number $k$, $4<k<5$, is to be determined.~Then
\[
{\small\text{$\frac 12\{f(h)+f(-h)\}=\frac 12\cdot (-1)^{m+n}\cdot (\pm |h|^k)$}}\quad\text{, for }{\small\text{$ h\in G$}},
\]
and {\small $\frac 12\{f(h)+f(-h)\}=0$} for {\small $h\notin G$}. Then {\small $f_{(4)}(0)=0$}, since {\small $\frac 12\{f(h)+f(-h)\}=o(h^4)$}, and {\small $f_{(6)}(0)$} does not exist, since {\small $\lim_{h\rightarrow 0} \frac 12\{f(h)+f(-h)\}$} is either $0$ or $\pm \infty $. We compute
{\small $\Delta_6^s(0,h;f)=f(3h)-6f(2h)+15f(h)-20f(0)+15f(-h)-6f(-2h)+f(-3h)
=(-1)^{m+n+1}3^k|h|^k-6\cdot (-1)^{m+n+1}2^k|h|^k+15\cdot (-1)^{m+n}|h|^k
=(-1)^{m+n+1}(3^k-6\cdot 2^k-15)|h|^k$}
and denote {\small $\varphi (k)=3^k-6\cdot 2^k-15$}. Since {\small $\varphi (4)\varphi (5)<0$}, by continuity, {\small $\varphi (k)=0$} for some $k$ in the open interval $(4,5)$. For that particular $k$, {\small $D_6^sf(0)=0$}, as we needed.

When $n=7$, let $G=\langle 3,5,7\rangle $ and  let $f$ be the function
\[
{\small\text{$f(x)=(-1)^{n+p}x^k$}},\quad\text{ for }{\small\text{$x=3^m5^n7^p\in G,$}}
\]
and {\small $f(x)=0$}, otherwise, where {\small $5<k<7$}. Then {\small $\Delta_7^s(0,2h;f)=f(7h)-7f(5h)+21f(3h)-35f(h)+\cdots $
$
=(-1)^{n+p+1}(7^k-7\cdot 5^k-21\cdot 3^k+35 )|h|^k$}, if {\small $h=3^m5^n7^p\in G.
$}
As usual, denote {\small $\varphi (k)=7^k-7\cdot 5^k-21\cdot 3^k+35$} and check that {\small $\varphi (5)\varphi (7)<0$}. The rest is the same as in the other two cases.

When $n=8$, let $G=\langle 2,3\rangle $ and let $f$ be the function
\[
{\small\text{$f(x)=(-1)^{m}x^k,$}}\quad\text{ for }{\small\text{$x=2^m3^n\in G,$}}
\]
and {\small $f(x)=0$}, for {\small $x\notin G$}, where {\small $7<k<8$}. Then {\small $\Delta_8^s(0,h;f)=f(4h)-8f(3h)+28f(2h)-56f(h)+70f(0)-\cdots$
$
=(-1)^{m}\left(4^k-8\cdot 3^k-28\cdot 2^k-56 \right)|h|^k\text{, for }h=2^m3^n\in G.
$}
The function {\small $\varphi (k)=4^k-8\cdot 3^k-28\cdot 2^k-56$} has {\small $\varphi (7)\varphi (8)<0$}, and the rest is folklore.
\end{proof}

Following the same method as in the proof of Theorem~\ref{T3.2} for $n=5,6,7,8$, when~$n=9$ we would start by letting $G=\langle 3,5,7\rangle $ and then look for an expression of $f(x)$ of~the~form
\[
{\small \text{$f(x)=(-1)^{am+bn+cp}x^k$}},\quad\text{ if }{\small\text{$x=3^m5^n7^p\in G$,}}
\]
and {\small $f(x)=0$}, otherwise, where {\small $7<k<9$} and {\small $a,b,c\in \{0,1\}$}. Then {\small $\Delta_9^s(0,2h;f)=f(9h)-9f(7h)+36f(5h)-84f(3h)+126f(h)-\cdots $} will be of the~form
{\small\[
(-1)^{am+bn+cp}\left(\pm 9^k\pm 9\cdot 7^k\pm36\cdot 5^k\pm 84\cdot 3^k\pm 126 \right)|h|^k.
\]}Taking {\small $\varphi(k)$} as the expression in the parenthesis, the only choices for {\small $\varphi(k)$} with the property that {\small $\varphi (7)\varphi (9)<0$} are {\small $\varphi(k)=\pm (9^k-9\cdot 7^k+36\cdot 5^k-84\cdot 3^k-126 )$}. Unfortunately, no choice for {\small $a,b,c\in \{0,1\}$} leads to either expression, and so the method in Theorem~\ref{T3.2} does not extend to the {\small $n=9$} case. In this way, part (ii) of Conjecture~B remains open for~{\small $n\geq 9$}.

\subsection*{Acknowledgment.} We are thankful to the anonymous reviewer for his/her careful reading, suggestions for improvement, and great insight.

\bibliographystyle{plain}

\begin{thebibliography}{ACCh1}
 \bibitem[AAC]{AAC} A. Ash, J. M. Ash and S. Catoiu, \textit{New definitions of continuity}, Real Anal. Exchange \textbf{40} (2014--15), no.~2, 403--420. 

\bibitem[As]{As} J. M. Ash, \textit{Generalizations of the Riemann derivative,} Trans. Amer. Math. Soc. \textbf{126} (1967), 181--199. 

\bibitem[As1]{As1} J. M. Ash, \textit{A characterization of the Peano derivative}, Trans. Amer. Math. Soc. \textbf{149} (1970), 489--501. 

\bibitem[As2]{As2} J. M. Ash, \textit{Remarks on various generalized
derivatives,} Special functions, partial differential equations, and harmonic analysis, pp. 25--39, 
Springer Proc. Math. Stat. \textbf{108}, Springer, Cham, 2014. 

\bibitem[AC]{AC}J. M. Ash and S. Catoiu, \textit{Quantum symmetric $L^p$ derivatives},
Trans. Amer. Math. Soc. \textbf{360} (2008), 959--987. 

\bibitem[AC1]{AC1} M. Ash and S. Catoiu, 
\textit{Multidimensional Riemann derivatives}, Studia Math. \textbf{235}  (2016),  no. 1, 87--100. 

\bibitem[AC2]{AC2} M. Ash and S. Catoiu, \textit{Characterizing Peano and symmetric derivatives and the GGR conjecture's solution}, Int. Math. Res. Notices, IMRN 2022, no. 10, 7893--7921.

\bibitem[ACCh]{ACCh} J. M. Ash, S. Catoiu and W. Chin, \textit{The classification of generalized Riemann derivatives,} Proc. Amer. Math. Soc. \textbf{146} (2018), no. 9, 3847--3862. 

\bibitem[ACCh1]{ACCH} J. M. Ash, S. Catoiu and W. Chin, \textit{The classification of complex generalized Riemann derivatives,} J.~Math. Anal. Appl. \textbf{502} (2021), no. 2, Article 125270. (40pp.) (doi.org/10.1016/j.jmaa.2021.125270)

\bibitem[ACCs]{ACCs} J. M. Ash, S. Catoiu and M. Cs\"{o}rnyei, \textit{Generalized vs. ordinary differentiation}, Proc. Amer. Math. Soc. \textbf{145} (2017), no. 4, 1553--1565. 

\bibitem[ACF]{ACF} J. M. Ash, S. Catoiu, and H. Fejzi\'{c}, \textit{Two pointwise characterizations of the Peano derivative}, preprint.

\bibitem[ACF1]{ACF1} J. M. Ash, S. Catoiu, and H. Fejzi\'{c}, \textit{A new proof of the GGR conjecture}, C. R. Math. Acad. Sci. Paris, to appear. 

\bibitem[ACR]{ACR}J.\ M. Ash, S. Catoiu, and R. R\'{\i}%
os-Collantes-de-Ter\'{a}n, \textit{On the nth quantum derivative,}\ J. Lond.
Math. Soc. \textbf{66} (2002), 114--130. 

\bibitem[AJ]{AJ} J. M. Ash and R. L. Jones, \textit{Mean value theorems for generalized Riemann derivatives}, Proc. Amer. Math. Soc. \textbf{101} (1987), no. 2, 263--271. 


\bibitem[C]{C} S. Catoiu, \textit{A differentiability criterion for continuous functions}, Monatsh. Math. \textbf{197} (2022), no. 2, 285--291. 

\bibitem[CF]{CF} S. Catoiu and H. Fejzi\'{c}, \textit{A generalization of the GGR conjecture}, Proc. Amer. Math. Soc., to appear.

\bibitem[dlVP]{dlVP} Ch. J. de la Vall\'ee Poussin, \textit{Sur l'approximation des fonctions d'une variable r\'eelle et de leurs d\'eriv\'ees par les p\^olynomes et les suites limit\'ees de Fourier}, Bull. Acad. Royale Belgique (1908), 193--254.

\bibitem[D]{D}A. Denjoy, \textit{Sur l'int\'egration des coefficients diff\'erentiels 
d'ordre sup\'erieur,} Fund. Math. \textbf{25} (1935), 273--326.

\bibitem[EW]{EW} M. J. Evans and C. E. Weil, \textit{Peano derivatives: A survey}, Real Anal. Exchange \textbf{7} (1981--82), no.1,~5-23. 

\bibitem[F]{F} H. Fejzi\'c, \textit{Decomposition of Peano derivatives}, Proc. Amer. Math. Soc. \textbf{119} (1993), no. 2, 599--609. 

\bibitem[F1]{F1} H. Fejzi\'c, \textit{Infinite approximate Peano derivatives}, Proc. Amer. Math. Soc. \textbf{131} (2003), no. 8, 2527--2536. 

\bibitem[FFR]{FFR} H. Fejzi\'c, C. Freiling and D. Rinne,
\textit{A mean value theorem for generalized Riemann derivatives},
Proc. Amer. Math. Soc. \textbf{136} (2008),  no. 2, 569--576.

\bibitem[FR]{FR} H. Fejzi\'c and D. Rinne, \textit{Peano path derivatives}, Proc. Amer. Math. Soc. \textbf{125} (1997), no. 9, 2651--2656.



\bibitem[GGR]{GGR} I. Ginchev, A. Guerraggio and  M. Rocca, \textit{Equivalence of (n+1)-th order Peano and usual derivatives for n-convex functions}, Real Anal. Exchange \textbf{25} (1999--00),  no. 2, 513--520.

\bibitem[GR]{GR} I. Ginchev, M. Rocca, \textit{On Peano and Riemann derivatives}, Rend. Circ. Mat. Palermo (2) \textbf{49}  (2000),  no.~3, 463--480.

\bibitem[HL]{HL1} P. D. Humke and M. Laczkovich, \textit{Convexity Theorems for Generalized Riemann Derivatives}, Real Anal. Exchange \textbf{15} (1989--90), no. 2, 652--674. 

\bibitem[HL1]{HL2} P. D. Humke and M. Laczkovich, \textit{Monotonicity theorems for generalized Riemann derivatives}, Rend. Circ. Mat. Palermo (2) \textbf{38} (1989), no. 3, 437--454. 

\bibitem[Ki]{Ki} A. Khintchine, \textit{Recherches sur la structure des fonctions mesurables,}
 Fund. Math. \textbf{9} (1927), 212--279.

\bibitem[LPW]{LPW} M. Laczkovich, D. Preiss and C. Weil, \textit{On unilateral and bilateral nth Peano derivatives}, Proc. Amer. Math. Soc. \textbf{99} (1987), no. 1, 129--134.

\bibitem[MZ]{MZ} J. Marcinkiewicz and A. Zygmund, \textit{On the differentiability of functions and summability of trigonometric series,} Fund. Math. \textbf{26} (1936), 1--43. 

\bibitem[MM]{MM} S. Mitra and S. N. Mukhopadhyay, \textit{Convexity conditions for generalized Riemann derivable functions}, Acta Math. Hungar.  \textbf{83}  (1999),  no. 4, 267--291.


\bibitem[P]{P} G. Peano, \textit{Sulla formula di Taylor}, Atti Acad. Sci. Torino \textbf{27} (1891--92), 40--46.


\bibitem[RAA]{RAA2} S. R\u{a}dulescu, P. Alexandrescu and D.-O. Alexandrescu, ,
\textit{The role of Riemann generalized derivative in the study of qualitative properties of functions},
Electron. J. Differential Equations 2013, no. 187,~14~pp. 

\bibitem[R]{R} B. Riemann, \textit{Ub\"er die Darstellbarkeit einer Funktion durch eine trigonometrische Reihe}, Ges. Werke, 2. Aufl., pp. 227--271. Leipzig, 1892.

\bibitem[SZ]{SZ} E. Stein and A. Zygmund, \textit{On the differentiability of functions}, Studia. Math. \textbf{23} (1964), 247--283.

\bibitem[T]{T1} B. S. Thomson, \textit{Monotonicity theorems}, Proc. Amer. Math. Soc. \textbf{83} (1981), 547--552. 

\bibitem[W]{W} C. E. Weil, \textit{Monotonicity, convexity and symmetric Peano derivatives}, Trans. Amer. Math. Soc. \textbf{231} (1976), 225--237.

\bibitem[Z]{Z} A. Zygmund, \textit{Trigonometric Series}, Vol. I, Cambridge University Press, 1959.
\end{thebibliography}

\end{document}